\documentclass[12pt]{amsart}

\usepackage{latexsym,amsmath,amsthm,amscd,amsfonts,diagrams,graphics}

\DeclareFontFamily{U}{rsf}{}
\DeclareFontShape{U}{rsf}{m}{n}{
  <5> <6> rsfs5 <7> <8> <9> rsfs7 <10-> rsfs10}{}
\DeclareMathAlphabet{\mathscr}{U}{rsf}{m}{n}

\DeclareMathAlphabet{\mathgth}{U}{euf}{m}{n}

\DeclareFontFamily{U}{cyr}{}
\DeclareFontShape{U}{cyr}{m}{n}{
  <5> wncyr5 <6> wncyr6 <7> wncyr7 <8> wncyr8 <9> wncyr9 <10-> wncyr10}{}
\DeclareMathAlphabet{\mathcyr}{U}{cyr}{m}{n}

\input cyracc.def

\newcommand{\Sha}{\mathcyr{\cyracc Sh}}

\newcommand{\cA}{{\mathscr A}}
\newcommand{\cB}{{\mathscr B}}

\newcommand{\cE}{{\mathscr E}}
\newcommand{\cF}{{\mathscr F}}
\newcommand{\cG}{{\mathscr G}}

\newcommand{\cI}{{\mathscr I}}

\newcommand{\cL}{{\mathscr L}}

\newcommand{\cN}{{\mathscr N}}
\newcommand{\cO}{{\mathscr O}}

\newcommand{\cU}{{\mathscr U}}

\newcommand{\bJ}{{\bar{J}}}

\newcommand{\vH}{\check{H}}

\newcommand{\dcE}{{\cE^{\cdot}}}

\newcommand{\FMYX}{\Phi_{Y\ra X}}

\newcommand{\sHom}{\underline{\mathrm{Hom}}}
\newcommand{\sExt}{\underline{\mathrm{Ext}}}
\newcommand{\sEnd}{\underline{\mathrm{End}}}

\newcommand{\D}{{\mathbf D}_{\mathrm{coh}}^b}

\newcommand{\chk}{{\scriptscriptstyle\vee}}
\newcommand{\R}{\mathbf{R}}
\newcommand{\Ld}{\mathbf{L}}
\newcommand{\lotimes}{\stackrel{\Ld}{\otimes}}

\DeclareMathOperator{\Spec}{Spec}

\DeclareMathOperator{\Aut}{Aut}

\DeclareMathOperator{\Hom}{Hom}

\DeclareMathOperator{\Pic}{Pic}

\DeclareMathOperator{\id}{id}

\DeclareMathOperator{\Cl}{Cl}
\DeclareMathOperator{\Br}{Br}

\DeclareMathOperator{\NS}{NS}

\DeclareMathOperator{\rk}{rk}
\DeclareMathOperator{\Ext}{Ext}

\newcommand{\ra}{\rightarrow}
\newcommand{\lra}{\longrightarrow}

\newcommand{\scdot}{\,\cdot\,}

\newcommand{\C}{\mathbf{C}}

\newcommand{\Q}{\mathbf{Q}}

\newcommand{\gMod}{\mathgth{Mod}}
\newcommand{\gCoh}{\mathgth{Coh}}

\newcommand{\iso}{\cong}
\newcommand{\Cech}{\v{C}ech }

\newcommand{\pj}{\mathbf{P}}

\newarrow{Equal} =====

\theoremstyle{plain}
\newtheorem{theorem}{Theorem}[section]
\newtheorem{lemma}[theorem]{Lemma}

\newtheorem{proposition}[theorem]{Proposition}
\newtheorem{conjecture}[theorem]{Conjecture}
\theoremstyle{definition}

\newtheorem{definition-theorem}[theorem]{Definition-Theorem}

\theoremstyle{remark}

\renewcommand{\phi}{\varphi}

\begin{document}

\author{Andrei C\u ald\u araru}

\title{Derived Categories of Twisted Sheaves on Elliptic Threefolds}

\date{}

\begin{abstract}
We construct an equivalence between the derived category of shea\-ves
on an elliptic threefold without a section and a derived category of
twisted sheaves (modules over an Azumaya algebra) on any small
resolution of its relative Jacobian.
\end{abstract}

\maketitle

\begin{section}{Introduction}
\label{sec:intro}

Equivalences of derived categories have attracted much interest in the
past few years, generated in part by applications to the study of
moduli spaces (\cite{MukK3}, \cite{MukAb}, \cite{BriEll}), and in part
by the conjectured relationship between derived categories and mirror
symmetry~(\cite{Kon}).  Recently there have been suggestions
(\cite[6.8]{Cal}, \cite{KapOrl}) that in order to obtain a good
description of physical phenomena, one should study not only the
derived category of usual sheaves but also the derived categories of
sheaves of modules over Azumaya algebras on the spaces involved.  From
a physical point of view, one needs $B$-fields in the construction of
conformal field theories, and these fields often have a component
which is an element of the Brauer group of the Calabi-Yau manifold
used.  Here we present an example of a Fourier-Mukai transform
involving a derived category of sheaves over an Azumaya algebra.  Few
such examples are known: a symplectic case was studied in~\cite{Pol},
and intersections of quadrics were studied in~\cite{Kap}, but none of
these relate to the Calabi-Yau property.

The simplest case of a Fourier-Mukai equivalence occurs between the
derived category of an elliptic curve $E$ and that of its dual
$\hat{E}$, and is induced by the Poincar\'e bundle on $E\times
\hat{E}$.  This has been generalized to a relative
situation~(\cite{BriEll}, \cite{BriEllK3}) as follows: starting with
an elliptic fibration $X\ra S$, one considers a relative moduli space
$Y\ra S$ of semistable sheaves on the fibers of $X\ra S$.  Under the
assumptions that $\dim X = \dim Y$ and that the moduli problem is
fine, one proves that $Y$ is smooth and that the universal sheaf
induces an equivalence of derived categories $\D(X) \iso \D(Y)$.

In this paper we study what happens when one removes the assumption
that the moduli problem is fine.  We work this out in the particular
case when $X\ra S$ is an elliptic fibration which does not possess a
section, and the moduli problem is the one that gives the {\em
relative Jacobian}, the moduli space of semistable sheaves of rank 1,
degree 0 along the fibers of $X\ra S$.  We analyze the ways in which
this moduli problem fails to be fine, and we obtain an equivalence of
derived categories (which involves sheaves of modules over an Azumaya
algebra) after making some necessary changes.

\subsection{}
\label{def:genellCY}
Let $f:X\ra S$ be an elliptic fibration, with $X$ and $S$ smooth
complex manifolds of dimensions 3 and 2, respectively, and satisfying
the following extra properties:
\begin{enumerate}
\item $f$ is flat (i.e.\ all fibers are 1-dimensional) and projective;
\item $f$ does not have any multiple fibers;
\item $f$ admits a multisection;
\item the discriminant locus $\Delta$ is a reduced, irreducible curve 
in $S$, having only nodes and cusps as singularities;
\item the fiber of $f$ over a general point of $\Delta$ is a rational 
curve with one node.
\end{enumerate}
We call such a map a {\em generic elliptic threefold}.  (For the
definition of elliptic fibration, multiple fibers, multisection,
discriminant locus, the reader is referred to~\cite[2.1]{DolGro}.)
The most interesting applications of our results are for Calabi-Yau
threefolds, but the Calabi-Yau condition is not needed for the results
in this paper to hold.

The reason for calling such an elliptic fibration ``generic'' is the
fact that in many families of elliptic threefolds the above properties
are shared by the general members of the family (especially when $X$
is Calabi-Yau, see~\ref{example:3sec} and~\ref{example:5sec}).  For
technical reasons, we shall only restrict our attention to generic
elliptic fibrations.

From here on fix a generic elliptic threefold $f:X\ra S$, and a
relatively ample line bundle $\cO_{X/S}(1)$ for $f$.

\subsection{}
\label{def:relJac}
Define the {\em relative Jacobian} $p:J\ra S$ of $f$ to be the
relative moduli space of semistable sheaves of rank 1, degree 0 on the
fibers of $f$.  (To be precise, rank 1, degree 0 is defined as having
the same Hilbert polynomial as the trivial line bundle.)  It is a
flat, projective elliptic fibration which has a natural section
$s:S\ra J$, obtained by associating to a point $t\in S$ the point
$[\cO_{X_t}]$ corresponding to the semistable sheaf $\cO_{X_t}$ on the
fiber $X_t$.  (It is not hard to see that $\cO_{X_t}$ is semistable
for all the fibers of a generic elliptic threefold,
Section~\ref{sec:dim1}.)

The above description of the relative Jacobian is based on Simpson's
construction~\cite{Sim} of relative moduli spaces of semistable
sheaves, a fact which allows us to analyze the geometry of $J$ by
studying moduli spaces of semistable sheaves on various degenerations
of the fibers of $f$.  These degenerations are well understood by work
of Miranda~\cite{Mir}, and moduli spaces of semistable sheaves on
these singular fibers have been studied in~\cite{OdaSes}.  Putting
together these results will give us a good understanding of the
geometry of $J$ (Section~\ref{sec:dim3}).

\subsection{}
In order to obtain an equivalence of derived categories between $X$
and $J$ one needs to find a good replacement for the notion of
universal sheaf.  This problem arises because in our situation a
universal sheaf does not exist on $X\times_S J$.  The main
contribution of this paper consists in the proposed solution to this
problem, as well as the consequences deduced from it.

There are two main obstructions to the existence problem: one is the
fact that there are properly semistable sheaves in the moduli problem
(on reducible fibers).  They are responsible for the apparition of
singularities in $J$ (a whole S-equivalence class of sheaves gets
contracted to a point).  This situation should be contrasted with
Bridgeland's result~\cite{BriEllK3} that when the moduli problem under
consideration is fine, the moduli space is smooth.  The approach we
take for solving this problem is to replace $J$ by an analytic small
resolution $\bJ$ of its singularities.  (In general, one can not
expect to be able to find an algebraic small
resolution,~\ref{subsec:noalgres}.  If one wanted to stay in the
algebraic realm, one would replace schemes by algebraic spaces, in the
sense of Artin.  We stick to the analytic situation for ease of
exposition.)

The other obstruction to the existence of a universal sheaf is the
fact that although one can find universal sheaves on $X\times_S U$ for
small enough open sets $U$ in $\bJ$, the lack of uniqueness of these
universal sheaves may prevent them from gluing together.  (More
precisely, these are pseudo-universal sheaves, parametrizing all the
stable sheaves and some of the semistable sheaves in the moduli
problem.)  In the particular situation under consideration this is
indeed the case, and the obstruction to this gluing is naturally an
element $\alpha$ of $\Br(\bJ)$, the Brauer group of $\bJ$.  We resolve
this problem replacing sheaves by $\alpha$-twisted sheaves in the
resulting equivalence of derived categories.  (A quick introduction to
derived categories of twisted sheaves is provided in
Section~\ref{sec:twsh}.)

\vspace{2mm}
\noindent
{\bf Theorem~\ref{thm:mainthm}.}  
{\em Let $X\ra S$ be a generic elliptic threefold, let $J\ra S$ be its
relative Jacobian, and let $\bJ\ra J$ be an analytic small resolution
of the singularities of $J$, with exceptional locus $E$.  Let
$\alpha\in\Br(\bJ)$ be the unique extension to $\bJ$ of the
obstruction to the existence of a universal sheaf on $X\times_S
(\bJ\setminus E)$.  Then there exists an $\alpha^{-1}$-twisted
pseudo-universal sheaf on $X\times_S \bJ$ whose extension by zero to
$X\times \bJ$ induces an equivalence of derived categories
\[ \D(X) \iso \D(\bJ, \alpha). \]
(By extending by zero we mean pushing forward by the natural inclusion
$X\times_S \bJ \hookrightarrow X\times \bJ$.  $\D(\bJ, \alpha)$
denotes the derived category of $\alpha$-twisted sheaves on $\bJ$.)
}

On a side note, it is worthwhile observing that there is a classical
construction which provides an element $\alpha\in \Br(\bJ)$ -- the
Ogg-Shafarevich theory of elliptic fibrations without a section.  It
is not hard to trace through the definitions and to see that the
element constructed by Ogg-Shafarevich theory agrees with the element
$\alpha$ we use in Theorem~\ref{thm:mainthm}.  Therefore the above
result can be seen as a generalization of Ogg-Shafarevich theory via
derived categories.
\vspace{1.5mm}

We note here a particularly striking consequence of
Theorem~\ref{thm:mainthm}:

\vspace{1.5mm}
\noindent
{\bf Theorem~\ref{thm:dadakCY}.}  
{\em Assume we are in the setup of Theorem~\ref{thm:mainthm}, and let
$n$ be the order of $\alpha$ in $\Br(\bJ)$.  Then we have
\[ \D(\bJ, \alpha) \iso \D(\bJ, \alpha^k), \]
for any $k$ coprime to $n$.
}

\vspace{1.5mm} 
The paper is structured as follows: in Section~\ref{sec:dim1} we give
a brief overview of moduli spaces of semistable sheaves on the
singular fibers that occur in generic elliptic threefolds.
Section~\ref{sec:dim3} is devoted to a geometric study of the relative
Jacobian and of its small resolutions.  The next section provides a
short introduction to the topic of twisted sheaves and their derived
categories (for complete details the reader should consult~\cite[Part
1]{Cal}).  In Section~\ref{sec:dereqv} we prove
Theorem~\ref{thm:mainthm}, and in a final section we discuss the
relationship of our result to Ogg-Shafarevich theory and
Theorem~\ref{thm:dadakCY}.

\vspace{1.5mm} 
\noindent
\textbf{Conventions.}  
We work over the field of complex numbers, and all the spaces
considered are analytic spaces.  The topology used is the analytic
topology, unless otherwise specified.  The same results hold in the
algebraic category, by replacing analytic spaces by algebraic spaces
and using the \'etale topology instead.

\vspace{1.5mm} 
\noindent
\textbf{Acknowledgments.}  
The results in this paper are part of my Ph.D.\ work, completed at
Cornell University.  I would like to thank my supervisor, Mark Gross,
for teaching me about twisted sheaves, elliptic fibrations, and
algebraic geometry in general, and for providing plenty of help and
encouragement.  The original idea of looking at this problem was his.

\end{section}

\begin{section}{The fiberwise picture}
\label{sec:dim1}

In this section we consider a generic elliptic threefold $f:X\ra S$,
and we review Miranda's results~\cite{Mir} about the kinds of
degenerate fibers of $f$ that can occur.  We then apply the results of
Oda and Seshadri~\cite{OdaSes}, which describe what the corresponding
degenerations are in the relative Jacobian.  For clarity, we also
sketch a proof of the fact that the moduli space of rank 1, degree 0
semistable sheaves on an $I_2$ curve is a nodal curve.

\subsection{}
\label{subsec:Mir}
The conditions that we have imposed on the elliptic fibration $f:X\ra
S$ highly restrict the possibilities of what singular curves can occur
as fibers of $f$.  The results in~\cite{Mir} imply that the fibers of
$f$ can be classified as follows:
\begin{enumerate}
\item[$(0)$] over $s\in S\setminus \Delta$, $X_s$ is a smooth elliptic
curve ($\Delta$ is the discriminant locus of $f$);
\item[$(I_1)$] over a smooth point $s$ of $\Delta$, $X_s$ is a rational curve
with one node;
\item[$(I_2)$] over a node $s$ of $\Delta$, $X_s$ is a reducible curve of type
$I_2$, i.e., two smooth $\pj^1$'s meeting transversely at two points;
\item[$(II)$] over a cusp $s$ of $\Delta$, $X_s$ is a rational curve with one 
cusp.
\end{enumerate}
Furthermore, in the case of the $I_2$ fiber, each component $C_i$ of the fiber
has normal bundle
\[ \cN_{C_i/X} \iso \cO_{C_i}(-1) \oplus \cO_{C_i}(-1). \]

\subsection{}
This classification allows us to do a case by case analysis of what
the moduli space of rank 1, degree 0 semistable sheaves on each one of
the four types of curves looks like.  We use Gieseker's notion of
stability, in its generalized form for pure sheaves used by Simpson
in~\cite{Sim}.  For a quick account of these notions,
see~\cite[Section 1.2]{BlueBook}.

Let $C$ be one of the four types of curves in the list
in~\ref{subsec:Mir}, arbitrarily polarized, and let $M$ be the moduli
space of semistable sheaves on $C$ whose Hilbert polynomial (with
respect to the polarization) equals that of $\cO_C$.  Then $M$ is
isomorphic to $C$ in cases $(0)$, $(I_1)$, $(II)$, while in case
$(I_2)$, $M$ is isomorphic to a rational curve with one node.  This is
an immediate consequence of the results in~\cite{OdaSes}, but for
clarity of the exposition we'll give a brief description of this
result, placing the emphasis on the geometry of case $(I_2)$.

\subsection{}
First recall the description of degree 0 line bundles on an elliptic
curve $C$.  All these line bundles are of the form $\cO_C(P-Q)$, with
$Q$ fixed and $P$ sweeping out $C$.  One can consider a universal
family $\cE$ on $C\times C$, whose restriction to $C\times \{P\}$ is
isomorphic to $\cO_C(P-Q)$ for any $P\in C$ (the second component of
the product).  One way to describe $\cE$ is
\[ \cE = \cO_{C\times C}(\Delta) \otimes \pi_1^*\cO_C(-Q), \]
where $\pi_1:C\times C\ra C$ is the projection onto the first factor
and $\Delta$ is the diagonal in $C\times C$.  Indeed, $\cE$ is
obviously flat over the second component, and its restriction
$\cE(P,Q) = \cE|_{C\times\{P\}}$ is isomorphic to 
\[ \cO(\Delta)|_{C\times\{P\}} \otimes \cO_C(-Q) = \cO_C(P)
\otimes \cO_C(-Q) = \cO_C(P-Q). \]

\subsection{}
In the above discussion we made use of the fact that $C$ was smooth,
and thus we were able to speak about $\cO_{C\times C}(\Delta)$.  We
want to generalize this to the other possible cases in~\ref{subsec:Mir}.

Let $C$ be any curve occurring in the classification~\ref{subsec:Mir},
let $Q$ be a fixed smooth point of $C$, and let $\Delta$ be the
diagonal in $C\times C$.  Define $\cO_{C\times C}(\Delta) =
\sHom(\cI_\Delta, \cO_{C\times C})$, the dual of the ideal sheaf of
$\Delta$.

Dualizing the exact sequence
\[ 0 \ra \cI_\Delta \ra \cO_{C\times C} \ra \cO_\Delta \ra 0, \]
one gets
\[ 0 \ra \cO_{C\times C} \ra \cO_{C\times C}(\Delta) \ra \sExt^1(\cO_\Delta,
\cO_{C\times C}) \ra 0. \]
A local computation shows that $\sExt^1(\cO_\Delta, \cO_{C\times C})$
is a line bundle on $\Delta$, and a deformation argument shows that
this line bundle is trivial.  We conclude that we have an exact
sequence
\[ 0 \ra \cO_{C\times C} \ra \cO_{C\times C}(\Delta)\ra \cO_\Delta \ra
0, \]
as expected.

Define $\cE = \cO_{C\times C}(\Delta) \otimes \pi_1^*\cO_C(-Q)$.  The
remainder of this section is devoted to the study of the properties of
$\cE$.

\subsection{}
\label{subsec:ext}
The restriction $\cE(P)$ of $\cE$ to a fiber $C\times \{P\}$ is the
unique non-trivial extension
\[ 0 \ra \cO_C(-Q) \ra \cE(P) \ra \cO_P\ra 0. \]
It is a torsion free sheaf for all $P\in C$, and for smooth $P$ it
equals $\cO_C(P-Q)$.  Thus it makes sense to consider the question of
stability for the sheaves $\cE(P)$.

This is the point where the analysis changes between the irreducible
fibers (cases $(0)$, $(I_1)$, $(II)$) and the reducible fibers (case
$(I_2)$).  In the first three situations, all the sheaves $\cE(P)$ are
stable.  However, in the $I_2$ case the stability of $\cE(P)$ depends
on the relative position of $P$ and $Q$.  This explains why in the
first three cases the moduli space is isomorphic to the curve itself,
while for the $I_2$ fibers the moduli space has a component
contracted.  Since these results are known (see~\cite{OdaSes}
or~\cite[6.3]{Cal}), from here on we only sketch what happens in the
$I_2$ case.

Let $C_1$ be the component of $C$ that $Q$ lies in, and let $C_2$ be
the other component.  If $P\in C_1\setminus C_2$, then $\cE(P)$ is
stable.  (This parallels the standard picture for degree 0 line
bundles on an elliptic curve.)  However, if $P\in C_2$, then
$\cE(P)$ is properly semistable, with composition factors
$\cO_{C_1}(-1)$ and $\cO_{C_2}(-1)$.  We can easily see how this
happens when $P$ is not a singular point of $C$, since then $\cE(P)$
is the line bundle $\cO_C(P-Q)$.  One has the exact sequence
\[ 0 \ra \cO_{C_2}(-S_1 - S_2) \ra \cO_C \ra \cO_{C_1} \ra 0, \]
where $S_1, S_2$ are the two singular points of $C$.  Tensoring it
with $\cO_C(P-Q)$ one obtains
\[ 0 \ra \cO_{C_2}(P-S_1-S_2) \ra \cO_C(P-Q) \ra \cO_{C_1}(-Q)  \ra
0, \] 
i.e.,\ 
\[ 0 \ra \cO_{C_2}(-1) \ra \cE(P) \ra \cO_{C_1}(-1) \ra 0 \]
because $C_1$ and $C_2$ are both isomorphic to $\pj^1$.  The reduced
Hilbert polynomials $p(\cO_{C_2}(-1);t)$ and $p(\cE(P);t)$ both equal
$t$, independent of the polarization of $C$.  (For a definition of the
reduced Hilbert polynomial, see~\cite[1.2.3]{BlueBook}.)  Hence
$\cO_{C_2}(-1)$ is a destabilizing subsheaf of $\cE(P)$, and the
composition factors are $\cO_{C_1}(-1)$ and $\cO_{C_2}(-1)$.  The
details of the same analysis for $P$ singular can be found
in~\cite[6.3.5]{Cal}.

\subsection{}
We conclude that for $P\in C_2$, all $\cE(P)$ are in the same
S-equivalence class.  If one denotes by $M$ the moduli space of rank
1, degree 0 semistable sheaves on $C$, then the family $\cE$ gives by
the universal property of $M$ a map $C\ra M$ which contracts the
component $C_2$ of $C$ to an image which is a rational curve with one
node (Figure~\ref{Fig1}).
\begin{figure}
\begin{center}
\includegraphics*{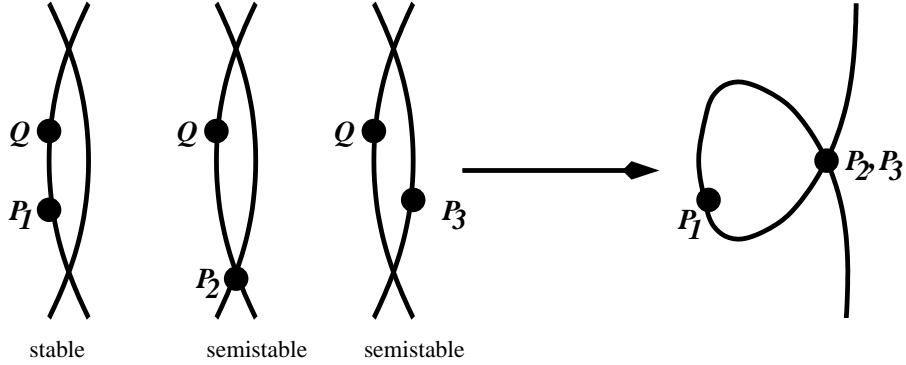}
\caption{The stability of the sheaves $\cE(P)$}
\label{Fig1}
\end{center}
\end{figure}
Since it is known that $M$ is a rational curve with one node, we
conclude that the map $C\ra M$ described above is onto.  

The way to think about this situation is that on $C\times M$ there is
no universal sheaf, because of the existence of properly semistable
sheaves.  There is not even a natural sheaf defined over all of
$C\times M$ and universal on $C\times M^{\mathrm{smooth}}$, because
there are many choices for the sheaf that would lie over the singular
point of $M$.  The solution is to replace $M$ by a ``blow-up'' of $M$
which is isomorphic to $C$, and which naturally parametrizes the
stable sheaves plus {\em some} semistable sheaves on $C$ (those of the
form $\cE(P)$ with $P\in C_2$).

\begin{proposition}
\label{prop:orthogonal}
If $C$ is of any one of the four types described in~\ref{subsec:Mir},
then the sheaves $\cE(P)$ satisfy
\[ \Ext^i(\cE(P_1), \cE(P_2)) = \left \{
\begin{array}{ll}
\C & \mbox{ if } P_1 = P_2 \mbox{ and } i=0 \\
0  & \mbox{ if } P_1 \neq P_2\mbox{ and any }i,
\end{array}
\right .
\]
for $P_1, P_2\in C$.  
\end{proposition}

In other words, in the sense of Fourier-Mukai transforms, the sheaves
in the family $\{\cE(P)\}_{P\in C}$ are ``almost'' mutually
orthogonal.  The ``almost'' refers to the fact that we may have
$\Ext^i(\cE(P), \cE(P)) \neq 0$ for $i>1$; this fact, however, is
purely a singular space phenomenon: we'll see later
(Proposition~\ref{prop:ceortho}) that if we embed $C$ in a smooth
space, the extensions by zero of the sheaves $\cE(P)$ are mutually
orthogonal in the proper sense.

This result shows that, from the point of view of derived categories,
using all the sheaves $\cE(P)$ (and not just the stable ones) is the
right thing to do.  Although we cannot speak of an equivalence of
derived categories induced by $\cE$ (because of the singularities of
$C$), morally $\cE$ should induce a Fourier-Mukai transform from $C$
to $C$.

\begin{proof}
Let $P$ be any point of $C$, and apply $\Hom(\scdot, \cO_C(-Q))$ to the
exact sequence
\[ 0 \ra \cO_C(-Q) \ra \cE(P) \ra \cO_P \ra 0,\]
to get 
\[ 0 \ra \Hom(\cE(P), \cO_C(-Q)) \ra \C \stackrel{\delta}{\lra} 
\Ext^1(\cO_P, \cO_C(-Q)). \] 
Since the extension is non-split, $\delta\neq 0$, and therefore
\[ \Hom(\cE(P), \cO_C(-Q)) = 0. \]
By Serre duality, $H^1(\cE(P)\otimes \cO_C(Q)) = 0$.  Since
\[ \chi(\cE(P)\otimes\cO_C(Q)) = 1 \]
we conclude that $\Hom(\cO_C(-Q), \cE(P)) = \C$.  We read this as
saying that, up to multiplication by a constant, there is a unique
morphism of $\cO_C(-Q)$ into $\cE(P)$.

Assume $P_1\neq P_2$, and let $f:\cE(P_1)\ra \cE(P_2)$ be any
homomorphism.  Compose $f$ with the inclusion $\cO_C(-Q)
\hookrightarrow \cE(P_1)$ to get a homomorphism $h:\cO_C(-Q)\ra
\cE(P_2)$.  Since there is a unique map (up to scalars) $\cO_C(-Q)\ra
\cE(P_2)$, $f$ can be restricted to a map $g$ that fits in the diagram
\[
\begin{diagram}[height=2em,width=2em]
0 & \rTo & \cO_C(-Q) & \rTo & \cE(P_1) & \rTo & \cO_{P_1} & \rTo & 0 \\
  &     & \dTo_g    &     & \dTo_f      &     & \dTo_0    &     &   \\
0 & \rTo & \cO_C(-Q) & \rTo & \cE(P_2) & \rTo & \cO_{P_2} & \rTo & 0.
\end{diagram}
\]
Since $\cO_C(-Q)$ is a line bundle, $g$ can be either an isomorphism
or zero.  If $g$ is an isomorphism, then the snake lemma shows that
$\ker f = \cO_{P_1}$, which is impossible because $\cE(P_1)$ is
torsion-free.  Therefore $g$ must be zero, and then $f$ is zero as
well.  We conclude that for $P_1\neq P_2$ we have
\[ \Hom(\cE(P_1), \cE(P_2)) = 0. \]

Note that $\cE(P_2)$ is locally free at $P_1$, so a local computation
shows that
\[ \Ext^i(\cO_{P_1}, \cE(P_2)) = \left \{ 
\begin{array}{ll}
\C & \mbox{ if } i = 1 \\
0 & \mbox{ otherwise.}
\end{array}
\right . \]
Applying $\Hom(\scdot, \cE(P_2))$ to the exact sequence
\[ 0 \ra \cO_C(-Q) \ra \cE(P_1) \ra \cO_{P_1} \ra 0 \]
gives 
\[ \Ext^i(\cE(P_1), \cE(P_2)) = \Ext^i(\cO_C(-Q), \cE(P_2)) \]
for all $i\geq 1$.  Since $C$ is Cohen-Macaulay with trivial dualizing
sheaf, and $\cO_C(-Q)$ is locally free, Serre duality shows that this
last group is 0 for $i\geq 2$, and is isomorphic to $\Hom(\cE(P_2),
\cO_C(-Q))^\chk$ for $i=1$.  Any non-trivial homomorphism $\cE(P_2)\ra
\cO_C(-Q)$ would give by composition a non-trivial one
$\cE(P_2)\ra\cE(P_1)$, contradicting our earlier result.  We conclude
that for $P_1\neq P_2$, $\Ext^i(\cE(P_1), \cE(P_2)) = 0$ for all $i$. 

Now let $P$ be any point of $C$, and apply $\Hom(\scdot, \cE(P))$ to
the exact sequence
\[ 0 \ra \cO_C(-Q) \ra \cE(P) \ra \cO_P \ra 0, \] 
to get
\[ 0 \ra \Hom(\cE(P), \cE(P))\ra \Hom(\cO_C(-Q), \cE(P)) = \C. \]
Since the first term is obviously non-trivial, we conclude that
\[ \Hom(\cE(P), \cE(P)) = \C. \]
\end{proof} 

\end{section}

\begin{section}{The relative Jacobian}
\label{sec:dim3}

In~\ref{def:relJac} we defined the relative Jacobian of a generic
elliptic threefold $X\ra S$ as the relative moduli space of semistable
sheaves on the fibers.  In this section we use the results in Section
2 to analyze the geometry of the relative Jacobian.

For completeness, we start by providing two examples of generic
elliptic threefolds without a section (in both these cases $X$ is also
Calabi-Yau).

\subsection{3-section}
\label{example:3sec}
Let $X$ be a general bidegree (3,3) hypersurface in $\pj^2\times
\pj^2$, considered with the projection $f:X\ra \pj^2$ to one of the
two factors of the product $\pj^2\times \pj^2$.  It is a Calabi-Yau
threefold, and the fibers of $f$ are degree 3 curves in $\pj^2$,
generically smooth, so that $f$ is an elliptic fibration.  The
discriminant locus is a reduced curve of degree 36, with 216 cusps and
189 nodes.  (This can be checked directly, using the software package
Macaulay~\cite{Mac}, or by using Euler characteristic calculations.)
Finally, it admits a 3-section (which can be taken to be a general
hyperplane section), and it does not admit a section (because, by the
Lefschetz theorem, $\Pic(X)$ is the restriction of $\Pic(\pj^2\times
\pj^2)$, and every divisor in $\pj^2\times\pj^2$ meets a fiber of $f$
with multiplicity divisible by 3).

\subsection{5-section}
\label{example:5sec}
In this example we construct a generic elliptic Calabi-Yau threefold
$X\ra \pj^2$, embedded in $\pj^2 \times \pj^4$.  Take coordinates
$x_0, \ldots, x_2$, $y_0, \ldots, y_4$ on $\pj^2\times \pj^4$, and let
$M$ be a generic $5\times 5$ skew-symmetric matrix whose $(i,j)$-th
entry is a polynomial of bidegree $(1,1)$ everywhere, except the last
row and column, where it is $(0,1)$.  According
to~\cite[0.1]{EisPopWal}, the $4\times 4$ Pfaffians of this matrix
define a degeneracy locus $X$, which has a symmetric locally free
resolution
\[ 0 \ra \cL \ra \cE \stackrel{\phi}{\lra} \cE^\chk(\cL) \ra 
\cO_{\pj^2\times\pj^4} \ra \cO_X \ra 0, \]
where
\begin{align*}
\cL & = \omega_{\pj^2\times\pj^4} = \cO(-3, -5), \\
\cE & = \bigoplus_{i=1}^5 \cO(a_i, -3), \\
(a_i) & = (-2, -2, -2, -2, -1),
\end{align*}
and the map $\phi$ is given by the matrix $M$.  Then it can be easily checked
using the results in~\cite{EisPopWal} that $X$ is a smooth Calabi-Yau
threefold.  

The projection of $X$ to $\pj^2$ is surjective and flat, and the fibers are
degree 5 curves in $\pj^4$ given by Pfaffians of a skew-symmetric $5\times 5$
matrix.  Therefore the fibers are (generically) elliptic curves, and it can be
checked by computer that this exhibits $X\ra \pj^2$ as a generic elliptic
fibration.  

The projection of $X$ to $\pj^4$ maps to a quintic threefold $Q$ in $\pj^4$,
contracting 52 lines and a conic, to 53 ordinary double points in $Q$.  It can
now be checked using standard techniques that the Picard number of $Q$ (and
therefore that of $X$) is 2.  Let $D$ and $H$ be pull-backs of hyperplane
sections from $\pj^2$ and $\pj^4$, respectively.  It is easy to compute
intersection numbers.  They are:
\[ D^3 = 0,\quad D^2H = 5,\quad DH^2 = 9,\quad H^3 = 5. \]
Since $D^2H$ and $DH^2$ are coprime, $D$ and $H$ must be primitive in $\NS(X)$,
so $\cO_X(D)$ and $\cO_X(H)$ generate $\Pic(X)$.  If $F=D^2$ is a fiber of
$X\ra \pj^2$, then we have $DF = 0$ and $HF=5$, so we conclude that $n=5$
(smallest degree of a multi-section).  One can take $H$ for a
multisection.

\subsection{}
\label{subsec:constrEi}
Let $f:X\ra S$ be a generic elliptic fibration, and let $J\ra S$ be
the relative Jacobian.  Since we have assumed that $f$ has no multiple
fibers, we can find a covering $\{U_i\}_{i\in I}$ of $S$ such that the
restriction $f_i:X_i=X\times_S U_i \ra U_i$ admits a section $s_i$ for
every $i\in I$.  Fix such a section $s_i$ for every $i\in I$.

On $X_i\times_S X_i$ we can consider the sheaf $\cE_i$ defined as
\[ \cE_i = \cI_\Delta^\chk \otimes \pi_1^* \cO_{X_i}(-s_i), \]
where $\Delta$ is the diagonal in $X_i\times_S X_i$, $\cI_\Delta^\chk$
is the dual 
\[ \cI_\Delta^\chk = \sHom_{X_i\times_S X_i}(\cI_\Delta,
\cO_{X_i\times_S X_i}), \]
and $\pi_1:X_i\times_S X_i \ra X_i$, is the projection onto the first
factor.  Then, as in Section~\ref{sec:dim1}, we find that $\cE_i$ is a
torsion-free sheaf, flat over the second component of the product
$X_i\times_S X_i$.

\subsection{}
For $P\in X_i$, let $X_P$ be the unique fiber of the elliptic
fibration in which $P$ lies, and let $Q$ be the point of $X_P$ that is
in the image of the section $s_i$.  Then the restriction $\cE_i(P) =
\cE_i|_{X_i\times_S \{P\}}$ is isomorphic to one of the sheaves
$\cE(P)$ considered in Section~\ref{sec:dim1}, and therefore it is
semistable on the fiber $X_P$.

This allows us to consider $X_i$ as a base space, parametrizing
sheaves on the fibers of $f$.  All the sheaves $\cE_i(P)$ are
semistable of rank 1, degree 0, so we get a natural map $\phi_i:X_i
\ra J_i = J\times_S U_i$, which only depends on the choice of the
section $s_i$.

Using our knowledge of the geometry of the fibers of $J$ we see that
the map $\phi_i$ is an isomorphism away from the $I_2$ fibers of $f$,
where it contracts the component of the fiber that is not hit by the
section $s_i$.  (One must remark here that since $s_i$ is a section,
it can not pass through any of the singular points of the singular
fibers of $X_i\ra U_i$.)  Since the contracted components have normal
bundle $\cO(-1)\oplus \cO(-1)$ by~\ref{subsec:Mir}, the contraction of
such a component gives an ordinary double point (ODP) in $J$.
Reversing the point of view, $X_i$ is a small resolution of the
singularities of $J_i$.

\subsection{}
\label{subsec:flop}
Each ODP naturally has {\em two} small resolutions (remember that we
are working in the analytic category) which are related by a flop.  As
long as $U_i$ is small enough so that $J_i$ contains at most one ODP
(which we can achieve since the $I_2$ fibers are isolated), both
resolutions of $J_i$ are isomorphic over $U_i$ to $X_i$.  Indeed, we
can replace the section $s_i$ by another section $s_i'$, which meets
the component of the $I_2$ fiber that was contracted by $\phi_i$.
Using $s_i'$ instead of $s_i$ gives a new map $\phi_i':X_i\ra J_i$,
which is easily seen to be the flop of $\phi_i$.

\subsection{}
The sheaves $\cE_i(P)$ remain mutually orthogonal in the global
setting, as shown by the next proposition:

\begin{proposition}
\label{prop:ceortho}
For $P\in X_i$, let $j_P:X_P\ra X$ be the inclusion of the fiber $X_P$
into $X$, and let $\cE_i^0(P) = j_{P,*} \cE_i(P)$.  Then the sheaves
$\cE_i^0(P)$ are mutually orthogonal, in the sense that
\[ \Ext^j_{X}(\cE_i^0(P_1), \cE_i^0(P_2)) = \left \{
\begin{array}{ll}
\C & \mbox{ if } P_1 = P_2 \mbox{ and } j=0, \\
0 & \mbox{ if } P_1\neq P_2 \mbox{ or } j>3.
\end{array}
\right .
\]
\end{proposition}

\begin{proof}
The result follows at once from Proposition~\ref{prop:orthogonal}, the
fact that the projective dimension of $X$ is $3$, and the following
lemma.
\end{proof}

\begin{lemma}
\label{lemma:fiberortho}
Let $f:X\ra S$ be a morphism of schemes or analytic spaces, with $S$
of the form $\Spec R$ for a regular local ring $R$.  If $s$ is the
closed point of $S$, let $i:X_s\ra X$ be the inclusion into $X$ of the
fiber $X_s$ over $s$, and let $\cF, \cG$ be sheaves on $X_s$.  If
$\Ext^j_{X_s}(\cF, \cG) = 0$ for all $j$ then $\Ext^j_X(i_* \cF, i_*
\cG) = 0$ for all $j$.
\end{lemma}

\begin{proof}
(An adaptation of the proof of~\cite[7.2]{BriEllK3}.)  We have 
\[ \R\Hom^\cdot_X(i_* \cF, i_* \cG) = \R\Hom^\cdot_{X_s}(\Ld i^* i_* \cF, \cG) 
\]
by the adjunction of $\Ld i^*$ and $i_*$.  Furthermore,
\[ \Ld i^* i_* \cF = \cF \lotimes_{X_s} \Ld i^* i_* \cO_{X_s} \]
by the projection formula.  Since $S$ is smooth at $s$, writing down the Koszul
resolution for $\cO_s$ and pulling back via $f$ we get a free resolution of
$\cO_{X_s}$ on $X$ which can be used to compute $\Ld i^* i_* \cF$.  This gives
\[ H^q(\Ld i^* i_* \cF) = \cF \otimes \bigwedge^q \cO_{X_s}^{\oplus m} \]
where $m=\dim S$.  Now the hypercohomology spectral sequence
\begin{eqnarray*}
E_2^{p,q} = \Ext^p_{X_s}(H^q(\Ld i^* i_* \cF), \cG) & \implies &
H^{p+q}(\R\Hom^\cdot_{X_s}(\Ld i^* i_* \cF, \cG)) \\
& = & H^{p+q}(\R\Hom^\cdot_X(i_* \cF, i_* \cG)) \\
& = & \Ext^{p+q}_X(i_*\cF, i_*\cG) 
\end{eqnarray*}
proves the result.  
\end{proof} 
\end{section}

\begin{section}{Twisted sheaves and derived categories}
\label{sec:twsh}

We sketch here the definition and main properties of twisted sheaves.
The reader unfamiliar with the subject is referred to~\cite[Chapters 1
and 2]{Cal} or~\cite{CalDTw}.

\subsection{}
Let $X$ be a scheme or analytic space, and let $\alpha\in \vH^2(X,
\cO_X^*)$ be represented by a \Cech 2-cocycle, given along a fixed
open cover $\{U_i\}_{i\in I}$ by sections
\[ \alpha_{ijk}\in \Gamma(U_i\cap U_j\cap U_k, \cO_X^*). \]
An $\alpha$-twisted sheaf $\cF$ (along the fixed cover) consists of a
pair 
\[ (\{\cF_i\}_{i\in I}, \{\phi_{ij}\}_{i,j\in I}), \]
where $\cF_i$ is a sheaf on $U_i$ for all $i\in I$ and
\[ \phi_{ij}: \cF_j|_{U_i\cap U_j} \ra \cF_i|_{U_i\cap U_j} \]
is an isomorphism for all $i,j\in I$, subject to the conditions:
\begin{enumerate}
\item $\phi_{ii} = \id$;
\item $\phi_{ij} = \phi_{ji}^{-1}$;
\item $\phi_{ij}\circ \phi_{jk}\circ \phi_{ki} = \alpha_{ijk} \cdot \id$.
\end{enumerate}

The class of twisted sheaves together with the obvious notion of
homomorphism is an abelian category, denoted by $\gMod(X, \alpha)$,
the category of $\alpha$-twisted sheaves.  If one requires all the
sheaves $\cF_i$ to be coherent, one obtains the category of coherent
$\alpha$-twisted sheaves, denoted by $\gCoh(X, \alpha)$.

This notation is consistent, since one can prove that these categories
are independent of the choice of the covering
$\{U_i\}$~(\cite[1.2.3]{Cal}) or of the particular cocycle
$\{\alpha_{ijk}\}$~(\cite[1.2.8]{Cal}) (all the resulting categories
are equivalent to one another).

\subsection{}
\label{subsec:functors}
For $\cF$ an $\alpha$-twisted sheaf, and $\cG$ an $\alpha'$-twisted
sheaf, one can define $\cF\otimes \cG$ (which is an
$\alpha\alpha'$-twisted sheaf), as well as $\sHom(\cF, \cG)$ (which is
$\alpha^{-1}\alpha'$-twisted), by gluing together the corresponding
sheaves.  If $f:Y\ra X$ is any morphism, $f^*\cF$ is an
$f^*\alpha$-twisted sheaf on $Y$.  Finally, if $\cF\in \gMod(Y,
f^*\alpha)$, one can define $f_* \cF$, which is $\alpha$-twisted on
$X$.  It is important to note here that one can not define arbitrary
push-forwards of twisted sheaves.

These operations satisfy all the usual relations (adjointness of $f_*$
and $f^*$, relations between $\sHom$ and $\otimes$, etc.)

The category $\gMod(X, \alpha)$ has enough injectives, and enough
$\cO_X$-flats~(\cite[2.1.1, 2.1.2]{Cal}).

\subsection{}
In the particular case when $\alpha$ can be represented by a sheaf
$\cA$ of Azumaya algebras over $X$ (in other words $\alpha\in\Br(X)$,
see~\cite{GrBr} and~\cite[Chapter IV]{Mil}) we can give a more
intrinsic description of the categories $\gMod(X, \alpha)$ and
$\gCoh(X, \alpha)$: they are equivalent to $\gMod(\cA)$ and
$\gCoh(\cA)$, the categories of sheaves of modules (respectively of
coherent sheaves of modules) over $\cA$.  This equivalence is obtained
by first remarking that there is a natural $\alpha$-twisted locally
free sheaf $\cE$ such that $\sEnd(\cE) \iso \cA$ (obtained by writing
locally $\cA\iso \sEnd(\cE_i)$ for some vector bundle $\cE_i$, and
gluing the $\cE_i$'s together into $\cE$) and thus the functors
\begin{eqnarray*}
F:\gMod(X, \alpha) \ra \gMod(\cA) & \quad & F(\scdot) =
\scdot\otimes_{\cO_X} \cE^\chk, \\
G:\gMod(\cA) \ra \gMod(X, \alpha) & \quad & G(\scdot) =
\scdot\otimes_{\cA} \cE
\end{eqnarray*}
are inverse to one another by standard Morita theory results.

\subsection{}
If $R$ is a commutative ring, $A$ and $B$ are Azumaya algebras over
$R$, then the $R$-linear categories $\gMod$-$A$ and $\gMod$-$B$ are
equivalent if and only if $[A] = [B]$ as elements of $\Br(R)$.  On the
other hand, if $R$ is a $k$-algebra for some subring $k$ of $R$, any
automorphism $\rho:R\ra R$ over $k$ induces an equivalence of
$k$-linear categories $\rho^*:\gMod$-$A \ra \gMod$-$\rho^* A$ for any
Azumaya algebra $A$ over $R$.  Since the action of $\Aut_k(R)$ on
$\Br(R)$ may be non-trivial, one sees at once that there may exist
Azumaya algebras $A$ and $B$ over $R$ which are not equal in the
Brauer group, but which are $k$-linearly Morita equivalent.  However,
the results in~\cite{Ros} suggest that any $k$-linear Morita
equivalence can be made into an $R$-linear one by pulling back by an
automorphism.  Therefore it makes sense to make the following
conjecture:

\begin{conjecture}
\label{conj:MorAzu}
Let $X$ be a scheme or complex analytic space, and let $\cA$, $\cB$ be
sheaves of Azumaya algebras over $X$.  Then the $\C$-linear categories
$\gCoh(\cA)$ and $\gCoh(\cB)$ are equivalent if and only if there
exists an automorphism $\rho$ of $X$ such that $\rho^*[\cA] = [\cB]$,
where $[\cA]$, $[\cB]$ denote the classes of $\cA$, $\cB$ in $\Br(X)$.
\end{conjecture}

In other words, the conjecture claims that the set of Azumaya algebras
on $X$, modulo Morita equivalence, is precisely the quotient of
$\Br(X)$ by the action of the automorphism group of $X$.

\subsection{}
We are mainly interested in $\D(\gMod(X, \alpha))$, the derived
category of complexes of $\alpha$-twisted sheaves on $X$ with coherent
cohomology.  For brevity, we'll denoted it by $\D(X, \alpha)$.  Since
the category $\gCoh(X, \alpha)$ does not have locally free sheaves of
finite rank if $\alpha\not\in\Br(X)$, from here on we'll only consider
the case $\alpha\in\Br(X)$.

The technical details of the inner workings of $\D(X, \alpha)$ can be
found in~\cite{CalDTw} or~\cite[Chapter 2]{Cal}.  The important facts
are that one can define derived functors for all the functors
considered in~\ref{subsec:functors}, and they satisfy the same
relations as the untwisted ones (see for example~\cite[II.5]{HarRD}).
One can prove duality for a smooth morphism $f:X\ra Y$, which provides
a right adjoint
\[ f^!(\scdot) = \Ld f^*(\scdot) \otimes_X \omega_{X/Y}[n] \]
to $\R f_*(\scdot)$, as functors between $\D(Y, \alpha)$ and $\D(X,
f^* \alpha)$.  

\subsection{}
If $X$ and $Y$ are smooth schemes or analytic spaces,
$\alpha\in\Br(Y)$, and $\dcE\in\D(X\times Y, \pi_Y^*\alpha^{-1})$
(where $\pi_X$ and $\pi_Y$ are the projections from $X\times Y$ to $X$
and $Y$ respectively), we define the integral functor
\begin{align*}
\FMYX^\dcE & :\D(Y, \alpha) \ra \D(X), 
\intertext{given by}
\FMYX^\dcE(\scdot) & = \pi_{X,*}(\pi_Y^*(\scdot) \lotimes \dcE).
\end{align*}

The following criterion for determining when $\FMYX^{\dcE}$ is an
equivalence (whose proof can be found in~\cite{CalDTw}
or~\cite[3.2.1]{Cal}) is entirely similar to the corresponding ones
for untwisted derived categories due to Mukai~\cite{MukAb},
Bondal-Orlov~\cite{BonOrl} and Bridgeland~\cite{Bri}.

\begin{theorem}
\label{thm:equiv}
The functor $F=\FMYX^\dcE$ is fully faithful if and only if for each
point $y\in Y$,
\[ \Hom_{\D(X)}(F\cO_y, F\cO_y) = \C, \]
and for each pair of points $y_1, y_2\in Y$, and each integer $i$,
\[ \Ext^i_{\D(X)}(F\cO_{y_1}, F\cO_{y_2}) = 0 \]
unless $y_1=y_2$ and $0\leq i\leq \dim Y$.  (Here $\cO_y$ is the
skyscraper sheaf $\C$ on $y$, which is naturally an $\alpha$-sheaf.)

Assuming the above conditions satisfied, then $F$ is an equivalence of
categories if and only if for every point $y\in Y$,
\[ F\cO_y \lotimes \omega_X \iso F\cO_y.\]
\end{theorem}

\end{section}

\begin{section}{The derived equivalence}
\label{sec:dereqv}

\subsection{}
The goal of this section is to prove the main theorem:

\begin{theorem}
\label{thm:mainthm}
Let $X\ra S$ be a generic elliptic threefold, let $J\ra S$ be its
relative Jacobian, and let $\bJ\ra J$ be an analytic small resolution
of the singularities of $J$, with exceptional locus $E$.  Let
$\alpha\in\Br(\bJ)$ be the unique extension to $\bJ$ of the
obstruction to the existence of a universal sheaf on $X\times_S
(\bJ\setminus E)$.  Then there exists an $\alpha^{-1}$-twisted
pseudo-universal sheaf on $X\times_S \bJ$ whose extension by zero to
$X\times \bJ$ induces an equivalence of derived categories
\[ \D(X) \iso \D(\bJ, \alpha). \]
(By extending by zero we mean pushing forward by the natural inclusion
$X\times_S \bJ \hookrightarrow X\times \bJ$.)
\end{theorem}

\subsection{}
If $X$ is Calabi-Yau, then an immediate consequence of this theorem is
the fact that $K_\bJ = 0$.  Indeed, this follows from the uniqueness
of the Serre functor, which is constructed categorically and is thus
preserved by an equivalence of categories.

It is also worthwhile remarking that, through the conjectural
translation of physics into derived categories, one should interpret
this theorem as saying that the conformal field theory built on $X$
(with no discrete torsion in the $B$-field) is equivalent to the one
built on $\bJ$, with discrete torsion $\alpha$ turned on.

\subsection{}
\label{subsec:noalgres}
Let us start by explaining the statement of the theorem.  In
Section~\ref{sec:dim3} we saw that the relative Jacobian $J$ is
singular, having an ODP for each $I_2$ fiber in $X$, caused by the
contraction of a whole S-equivalence class of properly semistable
sheaves on $X$.  The local picture we have developed suggests that in
order to obtain a well-behaved replacement for the universal sheaf,
one should consider a small resolution of these ODP's.

Unfortunately, we can not expect such a small resolution to exist in
the algebraic category: in order to obtain one, we need to blow up a
Weil divisor which is not even $\Q$-Cartier, and therefore we need to
have $\rk \Cl(J) > \rk\Pic(J)$.  It is not hard to prove that $\rk
\Cl(J) = \rk \Cl(X) = \rk \Pic(X)$, and hence if $\rk \Pic(X) = 2$ (as
is the case in both examples~\ref{example:3sec}
and~\ref{example:5sec}) $J$ can not have an algebraic small resolution
(because $\Pic(J)$ has rank at least $2$ having a section and a
projective base).

\subsection{}
\label{subsec:pseudo-univ}
Therefore we are led to considering an analytic small resolution $\bJ$
of the singularities of $J$ (we can take any such resolution).  Since
the singular points of $J$ coincide with its properly semistable
points, we'll call the points in the exceptional locus of the map
$\bJ\ra J$ semistable as well; the other points of $\bJ$ (which are in
a 1-1 correspondence with the stable points of $J$) will be called
stable.  

Cover $S$ with open sets $U_i$ small enough that there is at most
one $I_2$ fiber in $X_i = X\times_S U_i$, and such that a section
$s_i:U_i\ra X_i$ exists.  Let $J_i=J\times_S U_i$ and $\bJ_i =
\bJ\times_S U_i$.  Then, by~\ref{subsec:flop}, we can find an
isomorphism (of spaces over $S$) $\phi_i:\bJ_i \ra X_i$.
In~\ref{subsec:constrEi} we constructed a sheaf $\cE_i$ on
$X_i\times_S X_i$ which was a good replacement for the universal
sheaf.  Pulling back $\cE_i$ via the isomorphism $\id\times_S\,
\phi_i:X_i\times_S \bJ_i \ra X_i\times_S X_i$ we obtain a sheaf
$\cU_i$ on $X_i\times_S \bJ_i$ for each $i$.

We'll call such a sheaf $\cU_i$ a local {\em pseudo-universal} sheaf
for the moduli problem under consideration, because it parametrizes
the stable sheaves (just like a good universal sheaf should) but over
the semistable points of $\bJ$ it parametrizes some of the semistable
sheaves in the corresponding S-equivalence class.

\subsection{}
\label{subsec:defalpha}
We want to show that the local pseudo-universal sheaves $\cU_i$ can be
put together into a twisted (global) pseudo-universal sheaf $\cU$ on
$X\times_S \bJ$, and to understand the meaning of the actual twisting.

It is a general fact that, for any moduli problem of semistable
sheaves on a space $X$, if the stable part of the moduli space is
$M^s$, there exists a unique $\alpha\in\Br(M^s)$ such that a
$\pi_M^*\alpha^{-1}$-twisted sheaf exists on $X\times M^s$ (see, for
example,~\cite[3.3.2 and 3.3.4]{Cal}, as well as~\cite[Appendix
2]{MukK3}).  An easy explanation of this fact goes as follows: cover
$M^s$ with open sets $U_i$ small enough to have a local universal
sheaf $\cU_i$ over $X\times U_i$, and such that $\Pic(U_i\cap U_j) =
0$.  Then, from the universal property of the $\cU_i$'s, one concludes
that the restrictions $\cU_i|_{U_i\cap U_j}$ and $\cU_j|_{U_i\cap
U_j}$ are isomorphic.  Choosing isomorphisms $\phi_{ij}$ between them,
and using the fact that stable sheaves are simple, one finds that
there exists a \Cech 2-cocycle $\{\alpha_{ijk}\}$ such that
\[ \phi_{ij} \circ \phi_{jk} \circ \phi_{ki} = \alpha_{ijk}^{-1} \id, \]
therefore $(\{\cU_i\}, \{\phi_{ij}\})$ is an $\alpha^{-1}$-twisted
sheaf.  (The inverse is taken for notational convenience.)  The
uniqueness of $\alpha$ as a cohomology class is an easy consequence of
the universality of the $\cU_i$'s, and the fact that
$\alpha\in\Br(M^s)$ follows from their flatness.  The class $\alpha$
is called the {\em obstruction to the existence of a universal sheaf}
for the moduli problem under consideration.

The sheaves $\cU_i$ constructed in~\ref{subsec:pseudo-univ} restrict
to universal sheaves along the stable part of $\bJ$, so there exists a
unique $\alpha\in\Br(\bJ^s)$ which makes the collection
$\{\cU_i|_{\bJ^s}\}$ into a $\pi_{\bJ}^*\alpha^{-1}$-twisted sheaf on
$\bJ^s$.  But the condition that $\{\cU_i\}$ be a twisted sheaf only
concerns restrictions of these sheaves to intersections $(X_i\times_S
\bJ_i) \cap (X_j\times_S \bJ_j)$, and these are all in the stable part
of $\bJ$.  We conclude that there is a unique extension of $\alpha$ to
$\vH^2(\bJ, \cO_\bJ^*)$ that makes $\{\cU_i\}$ into a
$\pi_\bJ^*\alpha^{-1}$-twisted sheaf $\cU$, and from the flatness of
$\cU_i$ over $\bJ_i$ we conclude that $\alpha\in\Br(\bJ)$
by~\cite[3.3.4]{Cal}.  (To extend $\alpha$ from $\bJ^s$ to $\bJ$ we
could also have used the standard purity theorem for cohomological
Brauer groups,~\cite[III, 6.2]{GrBr}.)

\subsection{}
\label{subsec:alphagerbe}
In our particular case there is a more down-to-earth description of
the twisting.  First, note that the open sets $\bJ_{ij} = \bJ_i\cap
\bJ_j$ are not small enough to have trivial Picard group, but we can
force the Picard group to be generated by sections of the map
$\bJ_{ij} \ra U_{ij}=U_i\cap U_j$.  On $X_{ij} = X_i\cap X_j$ we have
the line bundle $\cO_{X_{ij}}(s_i-s_j)$, and we can consider its
pull-back $\cL_{ij}$ to $\bJ_{ij}$ via the isomorphism $\phi_i$.  The
collection $\{\cL_{ij}\}$ satisfies the following properties:
\begin{enumerate}
\item $\cL_{ii}$ is trivial;
\item $\cL_{ij} \iso \cL_{ji}^{-1}$;
\item $\cL_{ijk} = \cL_{ij} \otimes \cL_{jk} \otimes \cL_{ki}$ is
trivial, but no canonical trivialization exists;
\item $\cL_{ijk} \otimes \cL_{ijl}^{-1} \otimes \cL_{ikl} \otimes
\cL_{jkl}^{-1}$ is canonically trivial.
\end{enumerate}
These are precisely the necessary properties to make the collection
$\{\cL_{ij}\}$ represent an element $\alpha\in \vH^2(\bJ,
\cO_{\bJ}^*)$, which is viewed as a gerbe, i.e., an object obtained by
gluing together trivial gerbes via ``transition functions'' which are
line bundles (see~\cite{Gir},~\cite{Hit} and~\cite{Cha} for details).
The fact that the collection $\{\cU_i\}$ can be made into an
$\alpha^{-1}$-twisted sheaf $\cU$ for this $\alpha$ now follows from
the definitions.

\subsection{}
Having $\alpha$ and $\cU$ we can apply the criterion of
Theorem~\ref{thm:equiv} to the extension by zero $\cU^0$ of $\cU$ to
$X\times \bJ$.  (One can push forward, since both the twisting on
$X\times_S \bJ$ and on $X\times \bJ$ is pulled back from $\bJ$.)
Note that the sheaves $F\cO_x$ for $x\in\bJ$ are precisely
$\cE_i^0(P)$ of Proposition~\ref{prop:ceortho}, and they have the
orthogonality property required.  Therefore the integral transform
induced by $\cU^0$ is full and faithful.  Finally, since the canonical
class of the fibers of $X$ is trivial, the last condition of
Theorem~\ref{thm:equiv} is vacuously true, and hence $\Phi_{\bJ\ra
X}^{\cU^0}$ is an equivalence of categories between $\D(\bJ, \alpha)$
and $\D(X)$, thus finishing the proof of Theorem~\ref{thm:mainthm}.

\end{section}

\begin{section}{Ogg-Shafarevich theory; consequences of the main theorem}
\label{sec:conseq}

\subsection{}
We want to discuss the connection of the analysis in
Section~\ref{sec:dereqv} to Ogg-Shafare\-vich theory.  Recall that, in
its simplest form, Ogg-Shafarevich theory starts with an elliptic
fibration $J_U\ra U$ which has a section and does not have any
singular fibers, and establishes a 1-1 correspondence between the
Tate-Shafarevich group $\Sha_U(J_U) = \Br(J_U)/\Br(U)$ and isomorphism
classes (as spaces over $U$) of elliptic fibrations $X_U\ra U$ whose
relative Jacobian is isomorphic to $J_U\ra U$.  (One needs to be a bit
careful here: since $J_U\ra U$ has at least one non-trivial
automorphism -- the negation along the fibers -- one needs to rigidify
the isomorphism between $J_U$ and the relative Jacobian of $X_U$
further; otherwise one could not distinguish between $\alpha$ and
$\alpha^{-1}$ as elements of $\Sha_U(J_U)$.)  For a general discussion
of Ogg-Shafarevich theory and the Tate-Shafarevich group, the reader
is referred to ~\cite{DolGro}.

\subsection{}
\label{subsec:broggsha}
Let $X\ra S$ be a generic elliptic threefold without a section.  By
looking at its restriction $X_U$ over $U=S\setminus\Delta$, we obtain
an element $\alpha'\in \Sha_U(J_U)$.  Standard results
(see~\cite{DolGro}) show that there is a natural inclusion
$\Sha_U(J_U) \subseteq \Br'(\bJ)$ (the cohomological Brauer group of
$\bJ$), and thus $\alpha'$ can be regarded as an element in
$\Br'(\bJ)$.  Tracing through the definition of $\alpha'$ given by
Ogg-Shafarevich theory we find that it coincides precisely with the
one we gave for $\alpha$ as a gerbe, in~\ref{subsec:alphagerbe}.  We
conclude that the obstruction $\alpha$ to the existence of a universal
sheaf on $X\times_S \bJ$ matches the element $\alpha'$ constructed by
Ogg-Shafarevich theory.

\subsection{}
This partly answers an old question of Michael Artin: although one can
construct $J$ and $\alpha$ from $X$ in a purely algebraic way (using
the \'etale topology instead of the analytic one), there is no purely
algebraic construction that would take us from $J$ and $\alpha$ to
$X$.  The above analysis shows that one can take $X$ to be the
``spectrum'' of the derived category $\D(\bJ, \alpha)$, i.e.\ a space
$X$ with $\D(X)\iso\D(\bJ, \alpha)$.  The unsatisfying aspect of this
description is the fact that $X$ is not uniquely determined (see
below), and it is not obvious how to describe its elliptic structure
just in terms of its derived category.

\subsection{}
An interesting question one can ask now is if we can have 
\[ \D(\bJ, \alpha) = \D(\bJ, \beta) \]
for distinct $\alpha,\beta\in \Br(\bJ)$.  The above analysis, combined
with Theorem~\ref{thm:mainthm} allows us to find examples of this
phenomenon.  

\begin{theorem}
\label{thm:dadakCY}
Assume we are in the setup of Theorem~\ref{thm:mainthm}, and let
$n$ be the order of $\alpha$ in $\Br(\bJ)$.  Then we have
\[ \D(\bJ, \alpha) \iso \D(\bJ, \alpha^k), \]
for any $k$ coprime to $n$.
\end{theorem}

\subsection{}
Theorem~\ref{thm:dadakCY} is very surprising, in view of the following
result which shows that a similar phenomenon can not occur in a
local situation:
\begin{theorem}
\label{thm:dadaklocring}
Let $R$ be a commutative local ring, $A$ and $B$ Azumaya algebras over
$R$.  Then $A$ is $R$-linearly Morita equivalent to $B$ if and only if
$A$ is derived Morita equivalent to $B$, i.e.\ $\D(\gMod(A)) \iso
\D(\gMod(B))$ as $R$-linear triangulated categories.
\end{theorem}

\begin{proof}
\cite{Yek}.
\end{proof}

A reasoning similar to the one that led to
Conjecture~\ref{conj:MorAzu} shows that one should expect that over a
local $k$-algebra $R$ one has $\D(A) \iso \D(B)$ as $\C$-linear
triangulated categories if and only if there exists an automorphism of
$R$ over $k$ that takes $[A]$ to $[B]$ as elements of the Brauer group
of $R$.  In the global setting, Theorem~\ref{thm:dadakCY} shows that
this is no longer true: in Example~\ref{example:5sec} there is
no automorphism $\rho$ of the relative Jacobian $\bJ$ which takes
$\alpha$ to $\alpha^2$.  Indeed, let $\bJ$ be a small resolution of
the relative Jacobian of the elliptic fibration $X$ studied
in~\ref{example:5sec}, and let $\alpha\in\Br(\bJ)$ be the element that
corresponds to $X$.  Since $X$ has a 5-section, $\alpha$ has order 5
in $\Br(\bJ)$ by standard Ogg-Shafarevich theory.  Assume that $\rho$
is an automorphism of $\bJ$ that takes $\alpha$ to $\alpha^2$.  Since
$\rk \Pic(\bJ) = 2$, there is a unique elliptic fibration structure on
$\bJ$, so $\rho$ must act along the fibers of $\bJ\ra S$.  Restricting
it to the generic fiber $\bJ_\xi$ of $\bJ\ra S$, we obtain an
automorphism $\rho_\xi$ of the elliptic curve $\bJ_\xi$.  If it fixes
the origin it must have order divisible by 5 (since the map
$\alpha\mapsto \alpha^2$ has order 5), which is impossible
by~\cite[IV.4.7]{HarAG}.  Thus $\rho$ must correspond to translation
by a non-zero section $s$.  This section must be torsion since
otherwise $\rk\Pic(\bJ) > 2$.  Let $H$ be a general line in $S$, and
let $\bJ_H$ be the pull-back of $\bJ$ to $H$.  The map $\bJ_H\ra H$
has no multiple or reducible fibers, it is algebraic (we stayed away
from the small resolutions), and is not isotrivial, so we can
apply~\cite[Proposition 5.3.4 (ii)]{CosDol} to conclude that $\bJ_H\ra
H$ has no torsion sections, contradicting the existence of $s$.  We
conclude that there is no automorphism of $\bJ$ taking $\alpha$ to
$\alpha^2$, and therefore Theorem~\ref{thm:dadakCY} is a purely global
phenomenon.

\begin{proof}[Proof of Theorem~\ref{thm:dadakCY}]
Let $X^k\ra S$ be the relative moduli space of semistable sheaves of
rank 1, degree $k$ on the fibers of $X\ra S$, computed with respect to
a polarization $\cO_{X/S}(1)$ of fiber degree $n$.  (For the existence
of such a polarization see~\cite[Section 1]{DolGro}.)  This moduli
problem is fine (this is apparently a known fact,~\cite[4.2]{BriEll};
or one can prove it directly,~\cite[6.6.2]{Cal}), and using the
results in~\cite{BriEllK3} one concludes that $X^k$ is smooth, and the
universal sheaf on $X\times_S X^k$ induces an equivalence of derived
categories
\[ \D(X) \iso \D(X^k). \]

It is easy to see that $X^k$ is again a generic elliptic threefold:
$X$ and $X^k$ are locally isomorphic over open sets in $S$, by an
analysis similar to that of Section~\ref{sec:dim3}, simplified by the
fact that there are no contractions caused by semistable sheaves.
Thus $X$ and $X^k$ have the same singular fibers, and their
discriminant loci are the same, therefore $X^k$ is a generic elliptic
threefold.

The Jacobian of $X^k$ can be identified in a natural way with $J$, and
therefore we get an isomorphism
\[ \D(X^k) \iso \D(\bJ, \beta), \]
where $\beta$ is the element of $\Br(\bJ)$ that corresponds to $X^k\ra
S$.  A computation similar to the one in~\ref{subsec:alphagerbe} shows
that the restrictions of $\beta$ and $\alpha^k$ along $\bJ_{S\setminus
\Delta}$ are the same, so by~\ref{subsec:broggsha} we must have
$\beta=\alpha^k$ in $\Br(\bJ)$.  Therefore, we have
\[ \D(\bJ, \alpha) \iso \D(X) \iso \D(X^k) \iso \D(\bJ, \alpha^k). \]
\end{proof}

\subsection{}
Note that if $X$ is Calabi-Yau, then $X^k$ is Calabi-Yau as well.
Indeed, the equivalence of categories $\D(X) \iso \D(X^k)$ implies
$K_{X^k} = 0$ by the uniqueness of the Serre functor.  Also, the
universal sheaf induces an isomorphism
\[ H^{2,0}(X)\oplus H^{4,0}(X) \iso H^{2,0}(X^k) \oplus H^{4,0}(X^k), \]
by the standard technique of Mukai, and therefore $H^{2,0}(X^k) = 0$,
and we conclude that $X^k$ is Calabi-Yau.

This fact can be used to construct counterexamples to the Torelli
problem for Calabi-Yau threefolds.  See~\cite[6.7]{Cal}
and~\cite{CalTor} for details.
 
\end{section}

\bigskip \noindent
\small\textsc{Department of Mathematics and Statistics, \\
University of Massachusetts, \\
Amherst, MA 01003-4515, USA} \\
{\em e-mail: }{\tt andreic@math.umass.edu}


\begin{thebibliography}{GP}

\bibitem{Mac} 
Bayer, D., Stillman, M., {\em Macaulay: A system for computation in
algebraic geometry and commutative algebra}, 1982-1994.  Source and
object code available for Unix and Macintosh computers. Contact the
authors, or download from {\tt \small math.harvard.edu} via anonymous
ftp.

\bibitem{BonOrl}
Bondal, A., Orlov, D., {Semiorthogonal decompositions for algebraic 
  varieties}, preprint, alg-geom/9506012

\bibitem{Bri}
Bridgeland, T., {Equivalences of triangulated categories and Fourier-Mukai 
  transforms},  Bull. London Math. Soc. 31 (1999), no. 1, 25--34, (also 
  preprint, {\tt alg-geom/9809114})

\bibitem{BriEll}
Bridgeland, T., {Fourier-Mukai transforms for elliptic
  surfaces}, J. reine angew. math. 498 (1998) 115-133 (also preprint, 
  {\tt alg-geom/9705002})

\bibitem{BriEllK3}
Bridgeland, T., Maciocia, A., {Fourier-Mukai transforms for K3 fibrations}, 
  preprint, {\tt alg-geom/ 9908022}

\bibitem{Cal}
C\u ald\u araru, A., {\em Derived Categories of Twisted Sheaves on Calabi-Yau 
  Manifolds}, Ph.D.\ thesis, Cornell University (2000), also available at
  \newline
  {\tt http://www.math.umass.edu/\~{}andreic/thesis/maincorn.pdf}

\bibitem{CalDTw}
C\u ald\u araru, A., {Derived categories of twisted sheaves}, in preparation

\bibitem{CalTor}
C\u ald\u araru, A., {Counterexamples to Torelli via Fourier-Mukai Transforms},
   in preparation

\bibitem{Cha} 
Chatterjee, D.\ S., {\em On the Construction of abelian gerbs}, Ph.D.\ thesis,
  Cambridge (1998)

\bibitem{CosDol}
Cossec, F., Dolgachev, I., {\em Enriques Surfaces I}, Birkh\"auser, Boston, 
  (1989)

\bibitem{DolGro}
Dolgachev, I., Gross, M., {Elliptic three-folds I: Ogg-Shafarevich theory}, 
  J.\ Alg.\ Geom.\ 3 (1994), 39-80 (also preprint, {\tt alg-geom/9210009})

\bibitem{EisPopWal}
Eisenbud, D., Popescu, S., Walter, C., {Enriques surfaces and other 
  non-Pfaffian subcanonical subschemes of codimension 3}, MSRI preprint 037,
  (1999)

\bibitem{Gir}
Giraud, J., {\em Cohomologie non-ab\'elienne}, Grundlehren Vol.\ 179, 
  Springer-Verlag (1971)

\bibitem{Gro}
Gross, M., {Finiteness theorems for elliptic Calabi-Yau threefolds}, Duke 
  Math.\ J., Vol.\ 74, No.\ 2 (1994), 271-299

\bibitem{GrBr}
Grothendieck, A., {Le groupe de Brauer I--III,} in {\em Dix Expos\'es sur la 
  Cohomologie des Sch\'emas}, North-Holland, Amsterdam (1968), 46-188

\bibitem{HarAG} 
Hartshorne, R., {\em Algebraic Geometry}, Graduate Texts in Mathematics 
  Vol.\ 52, Springer-Verlag (1977)

\bibitem{HarRD} 
Hartshorne, R., {\em Residues and Duality}, Lecture Notes in Mathematics 
  Vol.\ 20, Springer-Verlag (1966)

\bibitem{Hit}
Hitchin, N.\ J., {Lectures on Special Lagrangian Submanifolds}, Lectures given
  at the ICTP School on Differential Geometry, April 1999, preprint, 
  {\tt math.DG/9907034}

\bibitem{BlueBook}
Huybrechts, D., Lehn, M., {\em Geometry of Moduli Spaces of Sheaves}, 
  Aspects in Mathematics Vol. E31, Vieweg (1997)

\bibitem{Kap} 

Kapranov, M.\ M., {On the derived category and $K$-functor of coherent sheaves
  on intersections of quadrics}, Izv. Akad. Nauk SSSR Ser. Mat. 52 (1988), no.
  1,186--199 (translation in Math. USSR-Izv. 32 (1989), no. 1, 191--204)

\bibitem{KapOrl} 
Kapustin, A., Orlov, D., {Vertex algebras, mirror symmetry, and
  D-branes: the case of complex tori}, preprint, {\tt hep-th/0010293}

\bibitem{Kon}
Kontsevich, M., {\it Homological algebra of mirror symmetry,}
Proceedings of the 1994 International Congress of Mathematicians I,
Birk\"auser, Z\"urich, 1995, p.~120 (also preprint, {\tt alg-geom/9411018})

\bibitem{Mil}
Milne, J.\ S., {\em \'Etale Cohomology}, Princeton Mathematical Series 33, 
  Princeton University Press (1980)

\bibitem{Mir}
Miranda, R., {Smooth models for elliptic threefolds}, in {\em Birational 
Geometry of Degenerations}, Birkh\"auser, (1983), 85-133

\bibitem{MukAb}
Mukai, S., {Duality between $D(X)$ and $D(\hat{X})$ with its application 
  to Picard sheaves}, Nagoya Math. J., Vol. 81 (1981), 153-175

\bibitem{MukK3}
Mukai, S., {On the moduli space of bundles on K3 surfaces, I.}, in {\em Vector
  Bundles on Algebraic Varieties}, Oxford University Press (1987), 341-413

\bibitem{OdaSes}
Oda, T., Seshadri, C.\ S., {Compactifications of the generalized Jacobian 
  variety}, Trans. Amer. Math. Soc., 253 (1979), 1-90

\bibitem{Pol}
Polishchuk, A., {Symplectic biextensions and a generalization of the 
  Fourier-Mukai transform}, Math. Res. Lett. 3 (1996), no. 6, 813--828

\bibitem{Ros}
Rosenberg, A.\ L., {The spectrum of abelian categories and reconstruction 
  of schemes}, in {\em Algebraic and Geometric Methods in Ring Theory},
  Marcel Dekker, Inc., New York, (1998), 255-274

\bibitem{Sim}
Simpson, C.\ T., {Moduli of representations of the fundamental group of a 
  smooth projective variety, I}, Publ. Math. IHES, 79 (1994), 47-129

\bibitem{Yek}
Yekuteli, A., private communication

\end{thebibliography}
\end{document}